\documentclass[12pt,reqno,a4paper]{amsart}
\usepackage{mathrsfs}
\usepackage{latexsym,amsmath,amsfonts,amssymb,amsthm}
\theoremstyle{plain}
\newtheorem{Thm}{Theorem}
\newtheorem{Lem}{Lemma}
\newtheorem{Cor}{Corollary}

\theoremstyle{definition}
\newtheorem*{Ack}{Acknowledgment}
\theoremstyle{remark}
\newtheorem*{Rem}{Remark}

\def\1{{\bf 1}}

\def\jacob #1#2{\genfrac{(}{)}{}{}{#1}{#2}}

\def\pmod #1{\ ({{\rm mod}}\ #1)}
\def\mod #1{\ {\rm mod}\ #1}

\begin{document}

\title{Lehmer-type congruences for lacunary harmonic sums modulo $p^2$}
\author{Hao Pan}
\address{Department of Mathematics, Nanjing University,
Nanjing 210093, People's Republic of China}
\email{haopan79@yahoo.com.cn} \subjclass[2000]{Primary 11A07;
Secondary 11B65}
\thanks{The author was supported by the National
Natural Science Foundation of China (Grant No.
10771135).}\keywords{}
\date{}\maketitle

\begin{abstract}
In this paper, we establish some Lehmer-type congruences for lacunary
harmonic sums modulo $p^2$.
\end{abstract}

\section{Introduction}
\setcounter{equation}{0} \setcounter{Thm}{0} \setcounter{Lem}{0}
\setcounter{Cor}{0}

The well-known Wolstenholme's harmonic series congruence asserts
that
\begin{equation}
\label{wolstenholme} \sum_{k=1}^{p-1}\frac{1}{k}\equiv 0\pmod{p^2}
\end{equation}
for each prime $p\geq 5$. With help of (\ref{wolstenholme}), Wolstenholme \cite{Wolstenholme62} proved
that
$$\binom{mp}{np}\equiv\binom{m}{n}\pmod{p^3}.
$$
for any $m,n\geq 1$ and prime $p\geq 5$. In 1938, Lehmer \cite{Lehmer38} discovered
an interesting congruence as follows:
\begin{equation*}
\sum_{j=1}^{(p-1)/2}\frac{1}{j}\equiv-\frac{2^p-2}{p}+\frac{(2^{p-1}-1)^2}{p}\pmod{p^2}
\end{equation*}
for each prime $p\geq 3$.

Define
$$
{\mathcal H}_{r,m}(n)=\sum_{\substack{1\leq k\leq n\\ k\equiv
r\pmod{m}}}\frac{1}{k}.
$$
Clearly, with help of \ref{wolstenholme}, Lehmer's congruence can be rewritten as
\begin{equation}
\label{lehmer2}
{\mathcal H}_{p,2}(p-1)\equiv\frac{2^{p-1}-1}{p}-\frac{(2^{p-1}-1)^2}{2p}\pmod{p}.
\end{equation}
In fact, Lehmer also proved three another
congruences in the same flavor:
\begin{equation}
\label{lehmer3} {\mathcal H}_{p,3}(p-1)\equiv\frac{3^{p-1}-1}{2p}-\frac{(3^{p-1}-1)^2}{4p}\pmod{p^2},
\end{equation}
\begin{equation}
\label{lehmer4} {\mathcal H}_{p,4}(p-1)\equiv\frac{3(2^{p-1}-1)}{4p}-\frac{3(2^{p-1}-1)^2}{8p}\pmod{p^2}
\end{equation}
and
\begin{equation}
\label{lehmer6} {\mathcal H}_{p,6}(p-1)\equiv\frac{2^{p-1}-1}{3p}+\frac{3^{p-1}-1}{4p}-\frac{(2^{p-1}-1)^2}{6p}-\frac{(3^{p-1}-1)^2}{8p}\pmod{p^2},
\end{equation}
where $p\geq 5$ is a prime. The proofs of
(\ref{lehmer2}),(\ref{lehmer3}),(\ref{lehmer4}) and (\ref{lehmer6})
are based on the values of Bernoulli polynomial $B_{p(p-1)}(x)$ at
$x=1/2,1/3,1/4,1/6$.

However, no another congruence for ${\mathcal H}_{p,m}(p-1)$ modulo $p^2$ is known,
partly since very few is known on the values of $B_{p(p-1)}(n/m)$ when $m\not=1,2,3,4,6$.
Some Lehmer-type congruences modulo $p$ (not modulo $p^2$!) have be proved in \cite{Williams82,Sun92,Sun93,Sun08,SunSun92,Sun02}.
In this paper, we shall investigate the Lehmer-type congruences modulo $p^2$.

Define
$$
{\mathcal T}_{r,m}(n)=\sum_{\substack{0\leq k\leq n\\ k\equiv
r\pmod{m}}}\binom{n}{k} \qquad\text{and}\qquad {\mathcal
T}_{r,m}^*(n)=\sum_{\substack{0\leq k\leq n\\ k\equiv
r\pmod{m}}}(-1)^k\binom{n}{k}.
$$
Clearly ${\mathcal T}_{r,m}^*(n)=(-1)^n{\mathcal T}_{n-r,m}^*(n)$
and
$$
{\mathcal T}_{r,m}^*(n)=\begin{cases} (-1)^r{\mathcal
T}_{r,m}(n)&\qquad\text{if }m\text{ is even},\\
(-1)^r({\mathcal T}_{r,2m}(n)-{\mathcal
T}_{m+r,2m}(n))&\qquad\text{if }m\text{ is odd}.
\end{cases}
$$
As we shall see soon, it is not difficult to show that
$$
{\mathcal H}_{r,m}(p-1)\equiv-\frac{{\mathcal
T}_{r,m}^*(p)-\delta_{r,m}(p)}p\pmod{p},
$$
where
\begin{equation}
\label{delta}
\delta_{r,m}(p)=\begin{cases}
1&\text{if }r\equiv 0\pmod{m},\\
-1&\text{if }r\equiv p\pmod{m},\\
0&\qquad\text{otherwise}.
\end{cases}
\end{equation}
\begin{Thm}
\label{t1}
Let $m\geq 2$ be an integer and let $p>3$ be a
prime with $p\not=m$. Then
\begin{align}
\label{t1e} {\mathcal H}_{p,m}(p-1)\equiv&-\frac{2{\mathcal
T}_{p,m}^*(p)+2}p+\frac{{\mathcal T}_{p,m}^*(2p)+2}{4p}\pmod{p^2}.
\end{align}
\end{Thm}
Let us see how (\ref{lehmer2}) follows from Theorem \ref{t1}. Clearly we have ${\mathcal T}_{0,2}^*(n)=2^{n-1}$ and ${\mathcal
T}_{1,2}^*(n)=-2^{n-1}$. Hence in view of (\ref{t1e}), for any prime
$p\geq 5$,
\begin{align*}
 {\mathcal H}_{p,2}(p-1)\equiv&-\frac{2{\mathcal
T}_{p,2}^*(p)+2}p+\frac{{\mathcal T}_{p,2}^*(2p)+2}{4p}=
\frac{2^p-2}p-\frac{2^{2p-1}-2}{4p}\pmod{p^2}.
\end{align*}

In \cite{Sun02}, Sun had showed that ${\mathcal T}_{r,m}(n)$ can be expressed in terms of some
linearly recurrent sequences with orders not exceeding $\phi(m)/2$, where $\phi$ is the Euler totient function.
Thus in view of Theorem \ref{t1}, for each $m$, we always have a Lehmer-type congruence for ${\mathcal H}_{p,m}(p-1)$ modulo $p^2$, involving some linearly recurrent sequences.

However, as we shall see later, (\ref{t1e}) is not suitable to derive (\ref{lehmer3}), (\ref{lehmer4}) and (\ref{lehmer6}). So we need the following theorem.
\begin{Thm}
\label{t2}
Let $m\geq 2$ be an integer and let $p>3$ be a
prime with $p\not=m$. Then
\begin{equation}
\label{t2e1} {\mathcal H}_{p,m}(p-1)\equiv-\frac{{\mathcal
T}_{p,m}^*(2p)+2}{4p}-\frac
p2\sum_{\substack{1\leq r\leq m\\
2r\not\equiv p\pmod{m}}}{\mathcal H}_{r,m}(p-1)^2\pmod{p^2}.
\end{equation}
\end{Thm}
When $m=3$, we have ${\mathcal T}_{p,3}^*(2p)=-2\times 3^{p-1}$ (cf.
\cite[Theorem 1.9]{Sun92} and \cite[Theorem 3.2]{Sun02}). Thus by
(\ref{t2e1}), we get
$$
{\mathcal H}_{p,3}(p-1) \equiv-\frac{{\mathcal
T}_{p,3}^*(2p)+2}{4p}-p\bigg(\frac{{\mathcal
T}_{p,3}^*(2p)+2}{4p}\bigg)^2=\frac{3^{p-1}-1}{2p}-\frac{(3^{p-1}-1)^2}{4p}\pmod{p^2},
$$
since
$$
{\mathcal H}_{0,3}(p-1)\equiv-{\mathcal
H}_{p,3}(p-1)\equiv\frac{{\mathcal T}_{p,3}^*(2p)+2}{4p}\pmod{p}.
$$

Let us apply Theorem \ref{t2} to obtain more Lehmer's type congruences. The Fibonacci numbers
$F_0,F_1,F_2,\ldots$ are given by
$F_0=0$, $F_1=1$ and $F_{n}=F_{n-1}+F_{n-2}$ for $n\geq 2$.
It is well-known that
$F_{p}\equiv\jacob{5}{p}\pmod{p}$ and
$F_{p-\jacob{5}{p}}\equiv0\pmod{p}$
for prime $p\not=2,5$, where $\jacob{\cdot}{p}$ is the Legendre symbol.
Williams \cite{Williams82} proved that
$$
\frac{2}{5}\sum_{1\leq k\leq
4p/5-1}\frac{(-1)^k}{k}\equiv\frac{F_{p-\jacob{5}{p}}}{p}\pmod{p}
$$
for prime $p\not=2,5$. Subsequently Sun and Sun \cite[Corollary 3]{SunSun92} proved that
\begin{equation}
\label{SunSun}
{\mathcal H}_{2p,5}(p-1)\equiv-{\mathcal
H}_{-p,5}(p-1)\equiv-\frac{F_{p-\jacob{5}{p}}}{2p}\pmod{p}.
\end{equation}
We have a Lehmer-type congruences as follows.
\begin{Thm}
\label{t3} Suppose that $p>5$ is a prime. Then
\begin{equation}
\label{t3e}
{\mathcal
H}_{p,5}(p-1)\equiv\frac{5^{\frac{p-1}{2}}F_{p}-1}{p}-\frac{5^{p-1}F_{2p-\jacob{5}{p}}-1}{4p}\pmod{p^2}.
\end{equation}
\end{Thm}
The Pell numbers $P_0,P_1,P_2,\ldots$ are given by
$P_0=0$, $P_1=1$ and $P_{n}=2P_{n-1}+P_{n-2}$ for $n\geq 2$.
We know $P_{p}\equiv\jacob{2}{p}\pmod{p}$ and $P_{p-\jacob{2}{p}}\equiv0\pmod{p}$
for every odd prime $p$. In \cite{Sun93}, Sun proved that
$$
(-1)^{\frac{p-1}2}\sum_{1\leq k\leq (p+1)/4}\frac{(-1)^k}{2k-1}\equiv-\frac{1}{4}\sum_{k=1}^{\frac{p-1}2}\frac{2^k}{k}\equiv\frac{P_{p-\jacob{2}{p}}}{p}\pmod{p}
$$
for odd prime $p$. Similarly, we have a Lehmer-type
congruence involving Pell numbers.
\begin{Thm}
\label{t4} Suppose that $p>3$ is a prime. Then
\begin{equation}
\label{t4e}
{\mathcal
H}_{p,8}(p-1)\equiv\frac{2^{2p-4}+2^{p-3}+2^{\frac{p-3}{2}}P_{p}-1}{p}-\frac{2^{4p-6}+2^{2p-4}+2^{p-2}P_{2p-\jacob{2}{p}}-1}{4p}\pmod{p^2}.
\end{equation}
\end{Thm}

We shall prove Theorems \ref{t1} and \ref{t2} in Section 2. And the proofs of Theorems \ref{t3} and \ref{t4} will be given in Section 3.

\section{Proof Theorems \ref{t1} and \ref{t2} }
\setcounter{equation}{0} \setcounter{Thm}{0} \setcounter{Lem}{0}
\setcounter{Cor}{0}

\begin{Lem} Suppose that $p$ is a prime. Then
\begin{equation}
\label{l1e1} \frac{1}{p}\sum_{\substack{1\leq k\leq p-1\\ k\equiv
r\pmod{m}}}(-a)^k\binom{p}{k}\equiv-\sum_{\substack{1\leq k\leq p-1\\
k\equiv r\pmod{m}}}\frac{a^k}{k}+p\sum_{\substack{1\leq j<k\leq
p-1\\ k\equiv r\pmod{m}}}\frac{a^k}{jk}\pmod{p^2} \end{equation} and
\begin{align}
\label{l1e2} \frac{1}{2p}\sum_{\substack{1\leq k\leq 2p-1,\ k\not=p\\
k\equiv r\pmod{m}}}(-a)^k\binom{2p}{k}
\equiv&-\sum_{\substack{1\leq k\leq p-1\\
k\equiv r\pmod{m}}}\frac{a^k}{k}-\sum_{\substack{1\leq k\leq p-1\\
k\equiv
2p-r\pmod{m}}}\frac{a^{2p-k}}{k}\notag\\&+2p\sum_{\substack{1\leq
j<k\leq p-1\\ k\equiv r\pmod{m}}}\frac{a^k}{jk}+2p\sum_{\substack{1\leq j<k\leq p-1\\
k\equiv 2p-r\pmod{m}}}\frac{a^{2-k}}{jk}\pmod{p^2}.
\end{align}
\end{Lem}
\begin{proof}
\begin{align*}
\frac{1}{p}\sum_{\substack{1\leq k\leq p-1\\ k\equiv
r\pmod{m}}}(-a)^k\binom{p}{k}=& \sum_{\substack{1\leq k\leq p-1\\
k\equiv
r\pmod{m}}}\frac{(-a)^k}{k}\prod_{j=1}^{k-1}\bigg(\frac{p}{j}-1\bigg)\\
\equiv&-\sum_{\substack{1\leq k\leq p-1\\ k\equiv
r\pmod{m}}}\frac{a^k}{k}+\sum_{\substack{2\leq k\leq p-1\\
k\equiv
r\pmod{m}}}\frac{a^k}{k}\sum_{j=1}^{k-1}\frac{p}{j}\pmod{p^2}.
\end{align*}
Similarly,
\begin{align*}
&\frac{1}{2p}\sum_{\substack{1\leq k\leq 2p-1,\ k\not=p\\ k\equiv
r\pmod{m}}}(-a)^k\binom{2p}{k}\\
=&\sum_{\substack{1\leq k\leq p-1\\
k\equiv
r\pmod{m}}}\frac{(-a)^k}{k}\binom{2p-1}{k-1}+\sum_{\substack{1\leq k\leq p-1\\
k\equiv 2p-r\pmod{m}}}\frac{(-a)^{2p-k}}{2p-k}\binom{2p-1}{k}.
\end{align*}
We have
$$
\sum_{\substack{1\leq k\leq p-1\\
k\equiv
r\pmod{m}}}\frac{(-a)^k}{k}\binom{2p-1}{k-1}\equiv-\sum_{\substack{1\leq
k\leq p-1\\ k\equiv
r\pmod{m}}}\frac{a^k}{k}+2p\sum_{\substack{1\leq j<k\leq p-1\\
k\equiv r\pmod{m}}}\frac{a^k}{jk}\pmod{p^2}.
$$
And
\begin{align*}
&\sum_{\substack{1\leq k\leq p-1\\
k\equiv 2p-r\pmod{m}}}\frac{(-a)^{2p-k}}{2p-k}\binom{2p-1}{k}\\
\equiv&\sum_{\substack{1\leq k\leq p-1\\
k\equiv 2p-r\pmod{m}}}\frac{a^{2p-k}}{2p-k}-2p\sum_{\substack{1\leq k\leq p-1\\
k\equiv 2p-r\pmod{m}}}\frac{a^{2-k}}{2p-k}\sum_{j=1}^{k}\frac{1}{j}\\
\equiv&-\sum_{\substack{1\leq k\leq p-1\\
k\equiv 2p-r\pmod{m}}}\bigg(\frac{a^{2p-k}}{k}+2p\cdot\frac{a^{2-k}}{k^2}\bigg)+2p\sum_{\substack{1\leq j\leq k\leq p-1\\
k\equiv 2p-r\pmod{m}}}\frac{a^{2-k}}{jk}\pmod{p^2}.
\end{align*}
We are done.
\end{proof}

Define
$$
{\mathcal S}_{r,m}(n)=\sum_{\substack{2\leq k\leq n\\ k\equiv
r\pmod{m}}}\frac{1}{k}\sum_{j=1}^{k-1}\frac{1}j.
$$
Substituting $a=1$ in (\ref{l1e1}), we get
\begin{Cor} Suppose that $m\geq 2$. Then
\begin{equation}
\label{c1e1} {\mathcal H}_{r,m}(p-1)\equiv-\frac{{\mathcal
T}_{r,m}^*(p)-\delta_{r,m}(p)}p+p{\mathcal S}_{r,m}(p-1)\pmod{p^2},
\end{equation}
where $\delta_{r,m}(p)$ is same as the one defined in
(\ref{delta}).
In particular,
\begin{equation}
\label{c1e2} {\mathcal H}_{p,m}(p-1)\equiv-\frac{{\mathcal
T}_{p,m}^*(p)+1}p+p{\mathcal S}_{p,m}(p-1)\pmod{p^2}.
\end{equation}
\end{Cor}
Substituting $r=p, p+m/2$ and $a=1$ in (\ref{l1e2}) and noting that $\binom{2p}{p}\equiv 2\pmod{p^3}$, we have
\begin{Cor} Suppose that $m\geq 2$. Then
\begin{equation}
\label{c2e1} {\mathcal H}_{p,m}(p-1)\equiv-\frac{{\mathcal
T}_{p,m}^*(2p)+2}{4p}+2p{\mathcal S}_{p,m}(p-1)\pmod{p^2}.
\end{equation}
And if $m$ is even, then
\begin{equation}
\label{c2e2} {\mathcal H}_{p+m/2,m}(p-1)\equiv-\frac{{\mathcal
T}_{p+m/2,m}^*(2p)}{4p}+2p{\mathcal S}_{p+m/2,m}(p-1)\pmod{p^2}.
\end{equation}
\end{Cor}
Combining (\ref{c1e2}) and (\ref{c2e1}), we get
$$
p{\mathcal S}_{p,m}(p-1)\equiv -\frac{{\mathcal
T}_{p,m}^*(p)+1}p+\frac{{\mathcal
T}_{p,m}^*(2p)+2}{4p}\pmod{p^2},
$$
and Theorem \ref{t1} easily follows.
\begin{Lem}
\label{l2}
$$
\sum_{r=1}^{m}{\mathcal
T}_{r,m}^*(n){\mathcal
T}_{r+s,m}^*(n)=(-1)^n{\mathcal
T}_{n+s,m}^*(2n).
$$
\end{Lem}
\begin{proof}
Let $\zeta$ be a primitive $m$-th root of unity.
Clearly,
\begin{align*}
{\mathcal
T}_{r,m}^*(n)=&\frac1m\sum_{k=0}^{n}(-1)^k\binom{n}{k}\sum_{t=1}^{m}\zeta^{(k-r)t}=\frac1m\sum_{t=1}^{m}\zeta^{-rt}(1-\zeta^t)^n.
\end{align*}
Hence
\begin{align*}
\sum_{r=1}^{m}{\mathcal
T}_{r,m}^*(n){\mathcal
T}_{r+s,m}^*(n)=&\sum_{r=1}^m\frac1{m^2}\sum_{1\leq t_1,t_2\leq m}\zeta^{-r(t_1+t_2)-st_2}(1-\zeta^{t_1})^n(1-\zeta^{t_2})^n\\
=&\frac{(-1)^n}m\sum_{t=1}^{m}\zeta^{-(n+s)t}(1-\zeta^{t})^{2n}
=(-1)^n{\mathcal
T}_{n+s,m}^*(2n).
\end{align*}
\end{proof}

By Lemma \ref{l2}, we have
$$
{\mathcal
T}_{p,m}^*(2p)=-\sum_{r=1}^{m}{\mathcal
T}_{r,m}^*(p)^2.
$$
Hence
\begin{align}
{\mathcal S}_{p,m}(p-1)\equiv&\frac{{\mathcal
T}_{0,m}^*(p)-1}{2p^2}-\frac{{\mathcal
T}_{p,m}^*(p)+1}{2p^2}-\frac{\sum_{r=1}^{m}{\mathcal
T}_{r,m}^*(p)^2-2}{4p^2}\notag\\
=&-\sum_{\substack{1\leq r\leq m\\ r\not\equiv 0,p\mod{m}}}\frac{{\mathcal
T}_{r,m}^*(p)^2}{4p^2}-\frac{({\mathcal
T}_{0,m}^*(p)-1)^2}{4p^2}-\frac{({\mathcal
T}_{p,m}^*(p)+1)^2}{4p^2}\notag\\
\equiv&-\frac{1}{4}\sum_{r=1}^m{\mathcal H}_{r,m}(p-1)^2\pmod{p}.
\end{align}
And since
${\mathcal H}_{p-r,m}(p-1)\equiv-{\mathcal H}_{r,m}(p-1)\pmod{p}$, ${\mathcal H}_{r,m}(p-1)\equiv0\pmod{p}$ provided that $2r\equiv p\pmod{m}$. So we also have
\begin{equation}
{\mathcal S}_{p,m}(p-1)\equiv-\frac{1}{4}\sum_{\substack{1\leq r\leq m\\ 2r\not\equiv p\pmod{m}}}{\mathcal H}_{r,m}(p-1)^2\pmod{p}.
\end{equation}
Thus by (\ref{c2e1}), Theorem \ref{t2} is concluded.

\section{Fermat's Quotient and Pell's Quotient}
\setcounter{equation}{0} \setcounter{Thm}{0} \setcounter{Lem}{0}
\setcounter{Cor}{0}

Let $L_n$ be the Lucas numbers given by
$L_0=2$, $L_1=1$ and  $L_{n}=L_{n-1}+L_{n-2}$ for $n\geq 2$.
We require the following result
of Sun and Sun on ${\mathcal T}_{r,10}(n)$.
\begin{Lem}{{\cite[Theorem 1]{SunSun92}}} \label{SunSunFib} Let $n$ be a positive odd integer. If $n\equiv
1\pmod{4}$, then
$$
\begin{array}{ccc}
&10{\mathcal
T}_{\frac{n-1}{2},10}(n)=2^n+L_{n+1}+5^{\frac{n+3}4}F_{\frac{n+1}2},\quad
&10{\mathcal T}_{\frac{n+3}{2},10}(n)=2^n-L_{n-1}+5^{\frac{n+3}4}F_{\frac{n-1}2},\\
&10{\mathcal
T}_{\frac{n+7}{2},10}(n)=2^n-L_{n-1}-5^{\frac{n+3}4}F_{\frac{n-1}2},\quad
&10{\mathcal
T}_{\frac{n+11}{2},10}(n)=2^n+L_{n+1}-5^{\frac{n+3}4}F_{\frac{n+1}2}.
\end{array}
$$
And if $n\equiv 3\pmod{4}$, then
$$
\begin{array}{ccc}
&10{\mathcal
T}_{\frac{n-1}{2},10}(n)=2^n+L_{n+1}+5^{\frac{n+1}4}L_{\frac{n+1}2},\quad
&10{\mathcal T}_{\frac{n+3}{2},10}(n)=2^n-L_{n-1}+5^{\frac{n+1}4}L_{\frac{n-1}2},\\
&10{\mathcal
T}_{\frac{n+7}{2},10}(n)=2^n-L_{n-1}-5^{\frac{n+1}4}L_{\frac{n-1}2},\quad
&10{\mathcal
T}_{\frac{n+11}{2},10}(n)=2^n+L_{n+1}-5^{\frac{n+1}4}L_{\frac{n+1}2}.
\end{array}
$$
Furthermore, for every odd $n$,
$$
10{\mathcal T}_{\frac{n+13}{2},10}(n)=2^n-2L_n.
$$
\end{Lem}
For each odd $n\geq 1$, since
$$
{\mathcal T}_{n,m}^*(2n)={\mathcal T}_{n,m}^*(2n-1)-{\mathcal
T}_{n-1,m}^*(2n-1)=-2{\mathcal T}_{n-1,m}^*(2n-1)
$$
and
$$
{\mathcal T}_{n+m,2m}^*(2n)={\mathcal T}_{n+m,2m}^*(2n-1)-{\mathcal
T}_{n+m-1,m}^*(2n-1)=-2{\mathcal T}_{n+m-1,m}^*(2n-1),
$$
by Lemma \ref{SunSunFib}, we get
\begin{equation}
{\mathcal T}_{n,5}^*(2n)=-2\cdot5^{\frac{n-1}2}F_{n}.
\end{equation}

Let $p>5$ be a prime. By (\ref{t2e1}),
$$
{\mathcal H}_{p,5}(p-1)\equiv-\frac{{\mathcal
T}_{p,5}^*(2p)+2}{4p}-p(H_{p,5}(p-1)^2+H_{2p,5}(p-1)^2)\pmod{p^2}.
$$
By (\ref{SunSun}), we have
\begin{align*}
{\mathcal
H}_{p,5}(p-1)\equiv&\frac{5^{\frac{p-1}2}F_{p}-1}{2p}-p\bigg(\bigg(\frac{F_{p-\jacob{5}{p}}}{2p}\bigg)^2+\bigg(\frac{5^{\frac{p-1}2}F_{p}-1}{2p}\bigg)^2\bigg)\\
\equiv&\frac{5^{\frac{p-1}2}F_{p}-1}{p}-\frac{5^{p-1}\big(F_{p-\jacob{5}{p}}^2+F_{p}^2\big)-1}{4p}\\
=&\frac{5^{\frac{p-1}2}F_{p}-1}{p}-\frac{5^{p-1}F_{2p-\jacob{5}{p}}-1}{4p}\pmod{p^2},
\end{align*}
where in the last step we use the fact $F_{2n-1}=F_n^2+F_{n-1}^2$.
Thus the proof of Theorem \ref{t3} is complete.
\begin{Rem} Similarly, we can get
\begin{align}
&\sum_{\substack{1\leq k\leq p-1\\ k\equiv p\pmod{5}}}\frac{(-1)^k}{k}\notag\\
\equiv&
\frac{5(2^{4p-1}-2^{2p+3})+12L_{4p}+L_{4p-4\jacob{5}{p}}-112L_{2p}-4L_{2p-2\jacob{5}{p}}+378}{400p}
\pmod{p^2}.\end{align}
\end{Rem}

Let $Q_n$ be the Pell-Lucas numbers given by
$Q_0=2$, $Q_1=2$ and $Q_{n}=2Q_{n-1}+Q_{n-2}$ for $n\geq 2$.
For ${\mathcal T}_{r,8}(n)$, Sun had proved that
\begin{Lem}{{\cite[Theorem 2.2]{Sun93}}}\label{SunPell} Let $n$ be a positive odd integer. If $n\equiv
1\pmod{4}$, then
$$
\begin{array}{ccc}
&8{\mathcal
T}_{\frac{n-1}{2},8}(n)=2^n+2^{\frac{n+1}2}+2^{\frac{n+7}4}P_{\frac{n+1}2},\quad
&8{\mathcal T}_{\frac{n+3}{2},8}(n)=2^n-2^{\frac{n+1}2}+2^{\frac{n+7}4}P_{\frac{n-1}2},\\
&8{\mathcal
T}_{\frac{n+7}{2},8}(n)=2^n-2^{\frac{n+1}2}-2^{\frac{n+7}4}P_{\frac{n-1}2},\quad
&8{\mathcal
T}_{\frac{n+11}{2},8}(n)=2^n+2^{\frac{n+1}2}-2^{\frac{n+7}4}P_{\frac{n+1}2}.
\end{array}
$$
And if $n\equiv 3\pmod{4}$, then
$$
\begin{array}{ccc}
&8{\mathcal
T}_{\frac{n-1}{2},8}(n)=2^n+2^{\frac{n+1}2}+2^{\frac{n+1}4}Q_{\frac{n+1}2},\quad
&8{\mathcal T}_{\frac{n+3}{2},8}(n)=2^n-2^{\frac{n+1}2}+2^{\frac{n+1}4}Q_{\frac{n-1}2},\\
&8{\mathcal
T}_{\frac{n+7}{2},8}(n)=2^n-2^{\frac{n+1}2}-2^{\frac{n+1}4}Q_{\frac{n-1}2},\quad
&8{\mathcal
T}_{\frac{n+11}{2},8}(n)=2^n+2^{\frac{n+1}2}-2^{\frac{n+1}4}Q_{\frac{n+1}2}.
\end{array}
$$
\end{Lem}
Thus we have
\begin{equation}
T_{n,8}^*(2n)=-2^{2n-3}-2^{n-2}-2^{\frac{n-1}2}P_{n}
\end{equation}
and
\begin{equation}
T_{n+4,8}^*(2n)=-2^{2n-3}-2^{n-2}+2^{\frac{n-1}2}P_{n}
\end{equation}
for odd $n\geq 1$.
Applying (\ref{t2e1}),
$$
{\mathcal H}_{p,8}(p-1)\equiv-\frac{{\mathcal
T}_{p,8}^*(2p)+2}{4p}-p\sum_{0\leq j\leq 3}{\mathcal
H}_{p+2j,8}(p)^2\pmod{p^2}.
$$
By (\ref{c2e1}) and (\ref{c2e2}), we have
$$
{\mathcal H}_{p,8}(p-1)\equiv-\frac{{\mathcal
T}_{p,8}^*(2p)+2}{4p}\pmod{p}
$$
and
$$
{\mathcal H}_{p+4,8}(p-1)\equiv-\frac{{\mathcal
T}_{p+4,8}^*(2p)}{4p}\pmod{p}.
$$
And in view of (\ref{c1e1}) and Lemma \ref{SunPell}.
\begin{align}
\label{hp26}
&{\mathcal H}_{p+2,8}(p-1)^2+{\mathcal
H}_{p+6,8}(p-1)^2\notag\\
\equiv&\begin{cases}
\sum_{i=0}^1p^{-2}\big(2^{p-3}-\jacob{2}{p}2^{\frac{p-5}2}+(-1)^i2^{\frac{p-5}4}P_{(p-\jacob{2}{p})/2}\big)^2\pmod{p}&\text{if
}p\equiv1\pmod{4},\\
\sum_{i=0}^1p^{-2}\big(2^{p-3}-\jacob{2}{p}2^{\frac{p-5}2}+(-1)^i2^{\frac{p-11}4}Q_{(p-\jacob{2}{p})/2}\big)^2\pmod{p}&\text{if
}p\equiv3\pmod{4}.
\end{cases}
\end{align}
\begin{Lem}
\begin{equation}
\label{pellcong}
\begin{array}{ccc}
P_{\frac{p-1}2}\equiv0\pmod{p},&P_{\frac{p+1}2}\equiv(-1)^{\frac{p-1}8}2^{\frac{p-1}4}\pmod{p},&\text{if }p\equiv 1\pmod{8},\\
P_{\frac{p-1}2}\equiv(-1)^{\frac{p-3}8}2^{\frac{p-3}4}\pmod{p},&P_{\frac{p+1}2}\equiv(-1)^{\frac{p+5}8}2^{\frac{p-3}4}\pmod{p},&\text{if }p\equiv 3\pmod{8},\\
P_{\frac{p-1}2}\equiv(-1)^{\frac{p-5}8}2^{\frac{p-1}4}\pmod{p},&P_{\frac{p+1}2}\equiv0\pmod{p},&\text{if }p\equiv 5\pmod{8},\\
P_{\frac{p-1}2}\equiv(-1)^{\frac{p+1}8}2^{\frac{p-3}4}\pmod{p},&P_{\frac{p+1}2}\equiv(-1)^{\frac{p+1}8}2^{\frac{p-3}4}\pmod{p},&\text{if
}p\equiv 7\pmod{8},
\end{array}
\end{equation} and
\begin{equation}
\label{pelllucascong}
\begin{array}{ccc}
Q_{\frac{p-1}2}\equiv(-1)^{\frac{p-1}8}2^{\frac{p+3}4}\pmod{p},&Q_{\frac{p+1}2}\equiv(-1)^{\frac{p-1}8}2^{\frac{p+3}4}\pmod{p},&\text{if }p\equiv 1\pmod{8},\\
Q_{\frac{p-1}2}\equiv(-1)^{\frac{p+5}8}2^{\frac{p+5}4}\pmod{p},&Q_{\frac{p+1}2}\equiv0\pmod{p},&\text{if }p\equiv 3\pmod{8},\\
Q_{\frac{p-1}2}\equiv(-1)^{\frac{p+3}8}2^{\frac{p+3}4}\pmod{p},&Q_{\frac{p+1}2}\equiv(-1)^{\frac{p-5}8}2^{\frac{p+3}4}\pmod{p},&\text{if }p\equiv 5\pmod{8},\\
Q_{\frac{p-1}2}\equiv0\pmod{p},&Q_{\frac{p+1}2}\equiv(-1)^{\frac{p+1}8}2^{\frac{p+1}4}\pmod{p},&\text{if
}p\equiv 7\pmod{8}.
\end{array}
\end{equation}
\end{Lem}
\begin{proof} The congruences in (\ref{pellcong}) were obtained by Sun \cite[Theorem
2.3]{Sun93}. And the congruences in (\ref{pelllucascong}) follows
from (\ref{pellcong}), by noting that $Q_{n}=2P_{n+1}-2P_{n}$ and
$Q_{n+1}=2P_{n+1}+2P_{n}$.
\end{proof}
Thus since $P_{(p-\jacob{2}{p})/2}Q_{(p-\jacob{2}{p})/2}=P_{p-\jacob{2}{p}}$, by (\ref{hp26}), we have
\begin{align*}
{\mathcal H}_{p+2,8}(p-1)^2+{\mathcal H}_{p+6,8}(p-1)^2
\equiv\frac{2^{p-1}(2^{\frac{p-1}2}-\jacob{2}{p})^2+P_{p-\jacob{2}{p}}^2}{8p}\pmod{p}.
\end{align*}
Observe that
$$
\frac{2^{p-1}-1}{p}=\frac{(2^{\frac{p-1}2}+\jacob2p)(2^{\frac{p-1}2}-\jacob2p)}{p}\equiv2\jacob2p\frac{2^{\frac{p-1}2}-\jacob2p}{p}\pmod{p}.
$$
Hence
\begin{align*}
{\mathcal H}_{p,8}(p-1)
\equiv&\frac{2^{2p-4}+2^{p-3}+2^{\frac{p-3}2}P_{p}-1}{p}-\frac{(2^{2p-3}+2^{p-2}+2^{\frac{p-1}2}P_{p})^2-4}{16p}\\
&-\frac{(2^{2p-3}+2^{p-2}-2^{\frac{p-1}2}P_{p})^2}{16p}-\frac{2^{p-1}(2^{\frac{p-1}2}-\jacob{2}{p})^2+2^{p-1}P_{p-\jacob{2}{p}}^2}{8p}\\
\equiv&\frac{2^{2p-4}+2^{p-3}+2^{\frac{p-3}2}P_{p}-1}{p}-\frac{2^{4p-6}+2^{2p-4}+2^{p-2}P_{2p-\jacob{2}{p}}-1}{4p}\pmod{p^2},
\end{align*}
by noting that $P_{p-\jacob{2}{p}}^2+P_{p}^2=P_{2p-\jacob{2}{p}}$.
This concludes the proof of Theorem \ref{t4}.\qed

\begin{Rem}
The Bernoulli polynomials $B_n(x)$ are given by
$$
\frac{te^{xt}}{e^{t}-1}=\sum_{n=0}^\infty\frac{B_n(x)}{n!}t^n.
$$
In particular, define the Bernoulli number $B_n=B_n(0)$. Granville and Sun \cite{GranvilleSun96} proved that
$$
B_{p-1}(\{p\}_5/5)-B_{p-1}\equiv\frac{5}{4p}F_{p-\jacob5p}+\frac{5^p-5}{4p}\pmod{p}
$$
and
$$
B_{p-1}(\{p\}_8/8)-B_{p-1}\equiv\frac{2}{p}P_{p-\jacob2p}+\frac{2^{p+1}-4}{p}\pmod{p}
$$
for prime $p\not=2,5$, where $\{p\}_m$ denotes the least non-negative residue of $p$ modulo
$m$. In \cite[Theorem 3.3]{Sun08}, Sun also proved that
$$
m{\mathcal H}_{p,m}(p-1)\equiv\frac{B_{2p-2}(\{p\}_m/m)-B_{2p-2}}{2p-2}-2\frac{B_{p-1}(\{p\}_m/m)-B_{p-1}}{p-1}\pmod{p^2}.
$$
Now using Theorems \ref{t3} and \ref{t4}
it is easy to deduce that
\begin{align}
\frac{B_{p(p-1)}(\{p\}_5/5)-B_{p(p-1)}}{5p(p-1)}
\equiv-\frac{5^{\frac{p-1}{2}}F_{p}-1}{p}+\frac{5^{p-1}F_{2p-\jacob{5}{p}}-1}{4p}\pmod{p^2}
\end{align}
for prime $p>5$, and
\begin{align}
&\frac{B_{p(p-1)}(\{p\}_8/8)-B_{p(p-1)}}{8p(p-1)}\notag\\
\equiv&-\frac{2^{2p-4}+2^{p-3}+2^{\frac{p-3}{2}}P_{p}-1}{p}+\frac{2^{4p-6}+2^{2p-4}+2^{p-2}P_{2p-\jacob{2}{p}}-1}{4p}\pmod{p^2}
\end{align}
for prime $p>3$.
\end{Rem}

\begin{Ack} The author is grateful to Professor Zhi-Wei Sun for his
helpful discussions on this paper.
\end{Ack}

\end{document}